\hsize 159.2mm
\vsize 246.2mm
\font\Bbb=msbm10

\font\bigrm=cmr17

\magnification=\magstep1
\def\X{\hbox{\Bbb X}}
\def\F{\hbox{\Bbb F}}
\def\R{\hbox{\Bbb R}}
\font\eightrm=cmr8

\footnote{}{ {\it AMS 2000 Subject Classifications.} 62G07, 62G20, 62F15.

 }

\centerline{\bigrm  Convergence Rates of Nonparametric }
\bigskip
\centerline{\bigrm   Posterior Distributions }
\vskip .5in
\centerline{\sl  Yang Xing    }
\vskip .2in
\centerline{\it  Swedish University of Agricultural Sciences}
\vskip .5in
\item{}{ \quad We study the asymptotic behavior of posterior distributions. We present general posterior convergence rate theorems, which extend several results on posterior convergence rates provided by Ghosal and Van der Vaart (2000), Shen and Wasserman (2001) and Walker, Lijor and Prunster (2007). Our main tools are the Hausdorff $\alpha$-entropy introduced by Xing and Ranneby (2008) and a new notion of prior concentration, which is a slight improvement of the usual prior concentration provided by Ghosal and Van der Vaart (2000). We apply our results to several statistical models.
}

\bigskip\bigskip
{\bf 1. Introduction.}\quad Recently, a major theoretical advance has occurred in the theory of Bayesian consistency for infinite-dimensional models. Schwartz (1965) first proved that, if the true density function is in the Kullback-Leibler support of the prior distribution, then the sequence of posterior distributions  accumulates in all weak neighborhoods of the true density function. It is known that the condition of positivity of prior mass on each Kullback-Leibler neighborhood in Schwartz's theorem is not a necessary condition. When one considers problems of density estimation, it is natural to ask for the strong consistency of Bayesian procedures. Sufficient conditions for the strong Hellinger consistency and for evaluating consistency rates have been currently developed by many authors. In this paper we study the problem of determining whether the posterior distributions accumulate in Hellinger neighborhoods of the true density function. The rate of this convergence can be measured by the size of the smallest shrinking Hellinger balls around the true density function on which posterior masses tend to zero as the sample size increases to infinity. By the fundamental works of Ghosal, Ghosh and Van der Vaart (2000) and Shen and Wasserman (2001), we know that the convergence rate of posterior distributions
is completely determined by two quantities: the rate of the metric entropy and the prior concentration rate. Roughly speaking, the rate of the metric entropy describes how large the model is, and the prior concentration rate depends on prior masses near the true distribution. Since the true distribution is unknown, the later assumption actually requires that the prior distribution spreads its mass "uniformly" over the whole density space. Another elegant approach for determination of the convergence rate was provided by Walker (2004), who obtained a sufficient condition for strong consistency by using summability of square root of prior probability instead of the metric entropy method. In this paper, in dealing with the rate of metric entropies we shall apply the Hausdorff $\alpha$-entropy introduced by Xing and Ranneby (2008), which is much smaller than  widely used metric entropies and the bracketing entropy. For some important prior distributions of statistical models the Hausdorff $\alpha$-entropies of all sieves are uniformly bounded, whereas it is generally impossible to get uniform boundedness of metric entropies of large sieves. The application of the Hausdorff $\alpha$-entropy leads refinements of several theorems on posterior convergence rates, for instance, the well known assumptions on metric entropies and summability of square root of prior probability have been weakened, which particularly yields that Theorem 5 of Ghosal et al.(2007b) is strengthened into Corollary 3 of this paper. To handle the prior concentration rate, we shall apply a new notion of prior concentration.
Our approach is a slight improvement of the  prior concentration provided by Ghosal et al.(2000), and moreover the proof of Lemma 1 on which the approach bases is quite simple. Finally, to get posterior convergence at the optimal rate $1/\sqrt n$, we give an extension of Ghosal et al.(2000, Theorem 2.4), in which the universal testing constant has been replaced by any fixed constant.

An outline of this paper is as follows. In Section 2 we define the Hausdorff $\alpha$-entropy with respect to a given prior and then present general theorems for the determination of posterior convergence rates. We also give a new approach to compute concentration rates. In Section 3 we apply our results to Bernstein polynomial priors, priors based on uniform distribution, log spline models and finite-dimensional models, which leads some improvements on known results for these models. The proofs of the main results are contained in Section 4.
\bigskip\bigskip
{\bf 2. Notations and Theorems.}\quad We consider a family of probability measures dominated by a $\sigma$-finite measure $\mu$ in $\X$, a Polish space endowed with a $\sigma$-algebra ${\cal X}$.
Let $X_1,X_2,\dots,X_n$ stand for an independent identically distributed (i.i.d.) sample of $n$ random variables, taking values in $\X$ and having a
common probability density function $f_0$ with respect to the measure  $\mu.$ Denote by $F^\infty_0$ the infinite product distribution of the probability distribution $F_0$ associated with $f_0$.
For two probability densities $f$ and $g$ we denote the Hellinger distance
$ H(f,g)= \Bigl(\int_{\X}\bigl(\sqrt{f(x)}-\sqrt{g(x)}\ \bigr)^2\mu(dx)\Bigr)^{1/2}$
and the Kullback-Leibler divergence
$K(f,g)=\int_{\X} f(x)\log {f(x)\over g(x)}\ \mu(dx).$
Assume that the space $\F$ of probability density functions is separable with respect to the Hellinger metric and that ${\cal F}$ is the Borel $\sigma$-algebra of $\F$.  Given a prior distribution $\Pi$ on ${\F}$, the posterior distribution $\Pi_n$ is a random probability measure with
 the following expression
 $$\Pi_n(A)=\Pi\bigl( A\,\big|\,X_1,X_2,\dots,X_n\bigr) ={\int_A\prod\limits_{i=1}^nf(X_i)\, \Pi(df)\over \int_{\F}\prod\limits_{i=1}^nf(X_i)\, \Pi(df)}={\int_A R_n(f)\, \Pi(df)\over \int_{\F} R_n(f)\, \Pi(df)}$$
for all measurable subsets $A\subset {\F}$, where $R_n(f)=\prod\limits_{i=1}^n\bigl\{f(X_i)/f_0(X_i)\bigr\}$ is the likelihood ratio. In other words, the posterior distribution $\Pi_n$  is the conditional distribution of $\Pi$ given the observations $X_1,X_2,\dots,X_n$. If the posterior distribution $\Pi_n$ concentrates on arbitrarily small neighborhoods of the true density function $f_0$ almost surely or in probability, then it is said to be consistent at $f_0$ almost surely and in probability respectively. Throughout this paper, almost sure convergence and convergence in probability should be understood as to be with respect to the infinite product distribution $F^\infty_0$ of $F_0$.

Our aim of this article is to present general theorems on  posterior convergence rates at $f_0$.
By the posterior convergence rate theorems of Ghosal, Ghosh and Van der Vaart (2000), we know that the prior concentration rate and the rate of metric entropy both completely determine the convergence rate of posterior distributions. More specifically, a key inequality to determine almost sure convergence rates of posterior distributions is
that for each $\varepsilon>0$,
$$\int_{\F} R_n(f)\, \Pi(df)\geq e^{-3\,n\,\varepsilon^2}\ \Pi\big( \,f: H(f_0,f)^2\,||f_0/f||_\infty<\varepsilon^2\,  \big)$$
almost surely for all sufficiently large $n,$ where $||g||_\infty$ stands for the supremum norm of the function $g$ on $\X$. This inequality was obtained by Ghosal et al.(2000, Lemma 8.4) under mild assumptions. It appears almost in all of papers handling strong convergence rates of posterior distributions. The reason is that in order to get the convergence rate of posterior distributions  one needs to find a suitable lower bound for the denominator in the expression of posterior distributions. This is successfully done in Ghosal et al.(2000), who suggested that the prior $\Pi$ puts sufficiently amount of mass around the true density function  $f_0$ in the sense: $\Pi\big( \,f: H(f_0,f)^2\,||f_0/f||_\infty<\tilde\varepsilon_n^2\,  \big)\geq e^{-n\,\tilde\varepsilon_n^2\,c}$ for some fixed constant $c$. Such a sequence $\{\tilde\varepsilon_n\}$ is referred to as the concentration rate of the prior $\Pi$ around $f_0$. Here we give a slightly stronger result. We introduce a modification of the Hellinger distance
$$ H_*(f_0,f)= \biggl(\int_{\X}\bigl(\sqrt{f_0(x)}-\sqrt{f(x)}\ \bigr)^2\Big({2\over 3}\,\sqrt{f_0(x)\over f(x)}+{1\over 3}\Big)\,\mu(dx)\biggr)^{1\over 2}.$$
It is clear that the inequality $||f_0/f||_\infty\geq 1$ holds for all density functions $f$ and $f_0$ such that the supremum is well-defined, and the quality holds if and only if $f=f_0$ almost surely. Observe also that $H_*(f_0,f)\not=H_*(f,f_0)$ and $3^{-1/2}\,H(f_0,f)\leq H_*(f_0,f)$. Moreover, we have
$$H_*(f_0,f)\leq  H(f_0,f)\,\big|\big|{2\over 3}\,\sqrt{f_0/ f}+{1\over 3}\big|\big|_\infty^{1/2}\leq  H(f_0,f)\,\big|\big|{f_0/ f}\big|\big|_\infty^{1/4}\leq H(f_0,f)\,\big|\big|{f_0/ f}\big|\big|_\infty^{1/2}$$
which yields
$$\bigl\{ f\in \F: H_*(f_0,f)\leq \tilde\varepsilon_n \bigr\}\supset \bigl\{ f\in \F: H(f_0,f)^2\,\sqrt{\bigl|\bigl| f_0\big/ f\bigr|\bigr|_\infty}\,<\tilde\varepsilon_n^2\ \bigr\}$$$$\supset \bigl\{ f\in \F: H(f_0,f)^2\,\bigl|\bigl| f_0\big/ f\bigr|\bigr|_\infty<\tilde\varepsilon_n^2\ \bigr\}.$$
The following simple lemma shows that, for $\tilde\varepsilon_n$ to be a
prior concentration rate, it is enough to assume $\Pi\bigl( W_{\tilde\varepsilon_n} \bigr)\geq e^{-n\,\tilde\varepsilon_n^2\,c_3},$ where 
$W_{\varepsilon}  = \bigl\{ f\in \F: H_*(f_0,f)\leq\varepsilon \bigr\}.$

\bigskip
\noindent {\bf Lemma 1.}  \it Let $\varepsilon>0$ and $c>0$. Then the inequality
$$F^\infty_0\Bigl(\ \int_{\F} R_n(f)\, \Pi(df)\leq e^{-n\,\varepsilon^2\,(3+2c)}\ \Pi\bigl( W_{\varepsilon}) \, \Bigr)\leq  e^{-n\,\varepsilon^2\, c} $$
holds for all $n$.
\rm
\bigskip
Lemma 1 provides a useful approach to compute prior concentration rates, particularly for models in which rate of convergence is governed by the prior concentration rate. It leads a simplification of the proof of Theorem 2.2 of Ghosal et al.(2000). Furthermore, we shall present  general posterior convergence rate theorems in which the well known assumptions on metric entropies and summability of square root of prior probability are also weakened.
We shall apply the Hausdorff $\alpha$-entropy $ J(\delta,{\cal G},\alpha)$ introduced by Xing et al.(2008). Denote by $L_\mu$  the space of all nonnegative integrable functions with the norm $||f||_1=\int_{\X} f(x)\,\mu(dx)$. Write $\log 0=-\infty$.
\bigskip
\noindent {\bf Definition.}   \it Let $\alpha\geq 0$ and ${\cal G}\subset {\F}$. For $\delta>0$, the Hausdorff $\alpha$-entropy $J(\delta,{\cal G},\alpha)$ with respect to the prior distribution $\Pi$ is defined as
$$J(\delta,{\cal G},\alpha)=\log\,\inf\ \sum\limits_{j=1}^N\,\Pi(B_j)^\alpha,$$
where the infimum is taken over all coverings $\{B_1,B_2,\dots,B_N\}$ of $\ {\cal G}$, where $N$ may take the value $\infty$, such that each $B_j$ is contained in some Hellinger ball $\{f:\,H(f_j,f)<\delta\}$ of radius $\delta$ and center at $f_j\in L_\mu$.
\rm
\bigskip
Note that the infimum can be equivalently taken over all partitions $\{P_1,P_2,\dots,P_N\}$ of $\ {\cal G}$ such that the Hellinger radius of each subset $P_j$ does not exceed $\delta$.
It was proved in Xing et al.(2008, Lemma 1) that the Hausdorff $\alpha$-entropy $J(\delta,{\cal G},\alpha)$ is an increasing subadditive function of ${\cal G}$ and satisfies $J(\delta,{\cal G},\alpha)\leq \log\,N(\delta,{\cal G})$ for all $\alpha\geq 0$, where $N(\delta,{\cal G})$ stands for the minimal number of Hellinger balls of radius $\delta$ needed to cover ${\cal G}.$ For $0\leq \alpha\leq 1$ and each $\cal G\subset \F$, we also obtained the following useful inequality
$$\Pi({\cal G})^\alpha\leq e^{J(\delta,{\cal G},\alpha)}\leq  \Pi({\cal G})^\alpha\,N(\delta,{\cal G})^{1-\alpha}.$$
Our first result in this paper is the following general theorem on posterior strong convergence rates. Denote
$A_\varepsilon=\bigl\{f:\,H(f_0,f)\geq \varepsilon\bigr\}.$
\bigskip
\noindent {\bf Theorem 1.}  \it Let $\{\bar\varepsilon_n\}_{n=1}^\infty$ and $\{\tilde\varepsilon_n\}_{n=1}^\infty$ be  positive sequences such that $n\,\min(\bar\varepsilon_n^2,\tilde\varepsilon_n^2)\to\infty$ as $n\to\infty$. Suppose that there exist constants  $c_1> 0,\,c_2>0,\,c_3\geq 0$, $0\leq \alpha<1$ and a sequence $\{{\cal G}_n\}_{n=1}^\infty$ of subsets on $\F$ such that $\sum\limits_{n=1}^\infty e^{-n\,\tilde\varepsilon_n^2\,c_2}<\infty$ and
$$\leqno \qquad (1)\quad \sum\limits_{n=1}^\infty  e^{J(\bar\varepsilon_n,{\cal G}_n,\alpha)-n\,\bar\varepsilon_n^2\,c_1}<\infty,$$
$$\leqno \qquad (2)\quad \sum\limits_{n=1}^\infty e^{n\,\tilde\varepsilon_n^2\,(3+3c_2+c_3)}\,\Pi(A_{\tilde\varepsilon_n}\setminus {\cal G}_n)<\infty ,$$
$$\leqno \qquad (3)\quad \Pi\bigl( f: H_*(f_0,f)\leq\tilde\varepsilon_n\ \bigr) \geq e^{-n\,\tilde\varepsilon_n^2\,c_3}.$$
Then for $\varepsilon_n=\max (\bar\varepsilon_n,\tilde\varepsilon_n)$ and each $r>2+\sqrt{2(3\alpha +2\alpha c_2+\alpha c_3+c_1)\over 1-\alpha}$, we have  $\Pi_n\bigl(A_{r\,\varepsilon_n}\bigr)\to 0$ almost surely as $n\to\infty$.
\rm
\bigskip
As direct applications we have
\bigskip
\noindent {\bf Corollary 1.}  \it Let $c_1\geq 0$, $c_2>0$, $c_3\geq 0$ and $0\leq\alpha<1$. Suppose that $\{\varepsilon_n\}_{n=1}^\infty$  is a positive sequence satisfying $\sum\limits_{n=1}^\infty e^{-n\,\varepsilon_n^2\,c_2}<\infty$ and suppose that there exists a sequence $\{{\cal G}_n\}_{n=1}^\infty$ of subsets on $\F$ such that
$$\leqno \qquad (1)\quad J(\varepsilon_n,{\cal G}_n,\alpha)\leq n\,\varepsilon_n^2\,c_1,$$
$$\leqno \qquad (2)\quad \Pi(\F\setminus {\cal G}_n)\leq e^{-n\,\varepsilon_n^2\,(3+3c_2+c_3)},$$
$$\leqno \qquad (3)\quad \Pi\bigl( f: H_*(f_0,f)\leq\varepsilon_n\ \bigr) \geq e^{-n\,\varepsilon_n^2\,c_3}.$$
Then for each $r>2+\sqrt{2(3\alpha +2\alpha c_2+\alpha c_3+c_1+ c_2)\over 1-\alpha}$, we have that $\Pi_n\bigl(A_{r\,\varepsilon_n}\bigr)\to 0$ almost surely as $n\to\infty$.
\rm
\bigskip
\noindent{\it Proof.} It is clear that all conditions of Theorem 1 are fulfilled if we let $\varepsilon_n=\bar\varepsilon_n=\tilde\varepsilon_n $ and replace the $c_1$ in Theorem 1 by $c_1+c_2$ of Corollary 1, and the proof is complete.
\bigskip
Corollary 1 extends Theorem 2.2 of Ghosal et al.(2000), in which they have stronger conditions than (1) and (3) of Corollary 1. It is probably worth mentioning that for several important prior distributions of infinite-dimensional statistical models, the quantities $J(\varepsilon_n,{\cal G}_n,\alpha)$ are uniformly bounded for all $n$ and hence condition (1) of Corollary 1 is trivially fulfilled, whereas general metric entropies of the sieve ${\cal G}_n$
grow to infinity as the sample size increases. A slightly different version of Theorem 1 is the following consequence, which is in fact an extension of Proposition 1 of Walker et all.(2007).
\bigskip
\noindent {\bf Corollary 2.}  \it Let $\{\varepsilon_n\}_{n=1}^\infty$  be a positive sequence such that $n\,\varepsilon_n^2\to\infty$ as $n\to\infty$.
Suppose that there exist constants  $c_1> 0$, $c_2> 0$, $c_3\geq 0$, $0\leq\alpha<1$  and a sequence $\{\cup_{j=1}^\infty {\cal G}_{nj}\}_{n=1}^\infty$ with ${\cal G}_{nj}\subset \F$ such that $\sum\limits_{n=1}^\infty e^{-n\,\varepsilon_n^2\,c_2}<\infty$ and
$$\leqno \qquad (1)\quad \sum\limits_{n=1}^\infty\,\sum\limits_{j=1}^\infty  N(\varepsilon_n, {\cal G}_{nj})^{1-\alpha}\,\Pi({\cal G}_{nj})^\alpha \,e^{-n\,\varepsilon_n^2\,c_1}<\infty;$$
$$\leqno \qquad (2)\quad \sum\limits_{n=1}^\infty e^{n\,\varepsilon_n^2\,(3+3c_2+c_3)}\,\Pi(A_{\varepsilon_n}\setminus \cup_{j=1}^\infty {\cal G}_{nj})<\infty ,$$
$$\leqno \qquad (3)\quad \Pi\bigl( f: H_*(f_0,f)\leq \varepsilon_n\ \bigr) \geq e^{-n\,\varepsilon_n^2\,c_3}.$$
Then for each $r>2+\sqrt{2(3\alpha +2\alpha c_2+\alpha c_3+c_1)\over 1-\alpha}$, we have that $\Pi_n\bigl(A_{r\,\varepsilon_n}\bigr)\to 0$ almost surely as $n\to\infty$.
\rm
\bigskip
\noindent{\it Proof.} Let ${\cal G}_n=\cup_{j=1}^\infty {\cal G}_{nj}$. We only need to verify condition (1) of Theorem 1 for such a sieve ${\cal G}_n$.
By Lemma 1 of Xing et al.(2008) we have
$$\sum\limits_{n=1}^\infty e^{J(\varepsilon_n,{\cal G}_n,\alpha)-n\,\varepsilon_n^2\,c_1}\leq \sum\limits_{n=1}^\infty\sum\limits_{j=1}^\infty
e^{J(\varepsilon_n,{\cal G}_{nj},\alpha)-n\,\varepsilon_n^2\,c_1}$$$$
\leq \sum\limits_{n=1}^\infty
\sum\limits_{j=1}^\infty N(\varepsilon_n, {\cal G}_{nj})^{1-\alpha}\,\Pi({\cal G}_{nj})^\alpha \,e^{-n\,\varepsilon_n^2\,c_1}<\infty.$$
Corollary 2 then follows from Theorem 1.
\bigskip
The assertion of Theorem 1 is an almost sure statement that the posterior distributions outside a Hellinger ball with a multiple of $\varepsilon_n$ as radius converge to zero almost surely.  Now we give an in-probability assertion under weaker conditions. Denote
$V(f,g)=\int_{\X} f(x)\bigl(\log {f(x)\over g(x)}\bigr)^2\ \mu(dx).$
\bigskip
\noindent {\bf Theorem 2.}  \it Let $\{\bar\varepsilon_n\}_{n=1}^\infty$ and $\{\tilde\varepsilon_n\}_{n=1}^\infty$ be positive sequences such that $n\,\min(\bar\varepsilon_n^2,\tilde\varepsilon_n^2)\to\infty$ as $n\to\infty$. Suppose that there exist constants $c_1> 0$, $c_2\geq 0$, $0\leq\alpha<1$ and a sequence $\{{\cal G}_n\}_{n=1}^\infty$ of subsets on $\F$ such that
$$\leqno \qquad (1)\quad J(\bar\varepsilon_n,{\cal G}_n,\alpha)-n\,\bar\varepsilon_n^2\,c_1\longrightarrow -\infty \qquad{\rm as}\quad n\to\infty,$$
$$\leqno \qquad (2)\quad e^{n\,\tilde\varepsilon_n^2\,(2+c_2)}\,\Pi(A_{\tilde\varepsilon_n}\setminus {\cal G}_n)\longrightarrow 0 \qquad{\rm as}\quad n\to\infty,$$
$$\leqno \qquad (3)\quad \Pi\bigl( f: K(f_0,f)<\tilde\varepsilon_n^2\ {\rm and}\ V(f_0,f)<\tilde\varepsilon_n^2\bigr)\geq e^{-n\,\tilde\varepsilon_n^2\,c_2}.$$
Then for $\varepsilon_n=\max (\bar\varepsilon_n,\tilde\varepsilon_n)$ and each $r>2+\sqrt{2(2\alpha +\alpha c_2+c_1)\over 1-\alpha}$, we have that $\Pi_n\bigl(A_{r\,\varepsilon_n}\bigr)\to 0$ in probability as $n\to\infty$.
\rm
\bigskip
Observe that the conditions (2) of Theorem 1 and Theorem 2 are only used to ensure that $\Pi_n(A_{\varepsilon_n}\setminus {\cal G}_n)\to 0$ as $n\to\infty$. So one can replace the conditions (2) of these theorems by  $\Pi_n(A_{\varepsilon_n}\setminus {\cal G}_n)\to 0$ as $n\to\infty$ almost surely and in probability respectively.
Theorem 2 is an extended version of Theorem 2.1 of Ghosal et al.(2001) and Theorem 1 of Walker et all.(2007). Furthermore, as a consequence of Theorem 2 we obtain the following slight improvement of Theorem 5 of Ghosal et al.(2007b).
\bigskip
\noindent {\bf Corollary 3.}  \it Let $\{\varepsilon_n\}_{n=1}^\infty$  be a positive sequence such that $n\,\varepsilon_n^2\to\infty$ as $n\to\infty$. Suppose that there exist constants $c_1> 0$, $c_2\geq 0$, $0\leq\alpha<1$  and a sequence $\{\cup_{j=1}^\infty {\cal G}_{nj}\}_{n=1}^\infty$ with ${\cal G}_{nj}\subset \F$. If
$$\leqno \qquad (1)\quad  e^{-n\varepsilon_n^2c_1}\,\sum\limits_{j=1}^\infty N(\varepsilon_n, {\cal G}_{nj})^{1-\alpha}\,\Pi({\cal G}_{nj})^\alpha \longrightarrow 0 \qquad{\rm as}\quad n\to\infty;$$
$$\leqno \qquad (2)\quad e^{n\,\varepsilon_n^2\,(2+c_2)}\,\Pi(A_{\varepsilon_n}\setminus \cup_{j=1}^\infty {\cal G}_{nj})\longrightarrow 0 \qquad{\rm as}\quad n\to\infty;$$
$$\leqno \qquad (3)\quad \Pi\bigl( f: K(f_0,f)<\varepsilon_n^2\ {\rm and}\ V(f_0,f)<\varepsilon_n^2\bigr) \geq e^{-n\,\varepsilon_n^2\,c_2},$$
then for each $r>2+\sqrt{2(2\alpha+\alpha c_2+c_1)\over 1-\alpha}$, we have that $\Pi_n\bigl(A_{r\,\varepsilon_n}\bigr)\to 0$ in probability as $n\to\infty$.
\rm
\bigskip
\noindent{\it Proof.} For ${\cal G}_n=\cup_{j=1}^\infty {\cal G}_{nj}$, by Lemma 1 of Xing et al.(2008) we get
$$e^{J(\varepsilon_n,{\cal G}_n,\alpha)-n\,\varepsilon_n^2\,c_1}\leq \sum\limits_{j=1}^\infty
e^{J(\varepsilon_n,{\cal G}_{nj},\alpha)-n\,\varepsilon_n^2\,c_1}
\leq e^{-n\,\varepsilon_n^2\,c_1}\,
\sum\limits_{j=1}^\infty N(\varepsilon_n, {\cal G}_{nj})^{1-\alpha}\,\Pi({\cal G}_{nj})^\alpha$$
which tends to zero as $n\to\infty$ and condition (1) of Theorem 2 holds. Then by Theorem 2 we conclude the proof.
\bigskip

The above theorems cannot yield a convergence rate $1/\sqrt n$ because of the assumption $n\,\varepsilon_n^2\to \infty$. Particularly, these theorems cannot well serve finite-dimensional models. Ghosal et al.(2000, 2007a) have obtained a nice theorem to handle such models. Denote $B_{\varepsilon_n^2}=\bigl\{ f: K(f_0,f)<\varepsilon_n^2\ {\rm and}\ V(f_0,f)<\varepsilon_n^2\bigr\}.$ Now our result is
\bigskip
\noindent {\bf Theorem 3.}  \it Let $\{\varepsilon_n\}_{n=1}^\infty$  be a positive sequence such that  $n\,\varepsilon_n^2$ are uniformly bounded away from zero, i.e., there exists a constant $c_0>0$ such that $n\,\varepsilon_n^2\geq c_0$ for all $n$  . Suppose that there exist constants  $0<\alpha<1$, $c_1< {1-\alpha\over 18}$ and a sequence $\{{\cal G}_n\}_{n=1}^\infty$ of subsets on $\F$ such that
$$\leqno \qquad (1)\quad {e^{2n\varepsilon_n^2}\, \Pi(A_{\varepsilon_n}\setminus {\cal G}_n)\over \Pi(B_{\varepsilon_n^2})}\longrightarrow 0 \qquad{\rm as}\quad n\to\infty,$$
$$\leqno \qquad (2)\quad \exp\Bigl( J\bigl({j\varepsilon_n\over 3},\bigl\{f\in {\cal G}_n:\,j\varepsilon_n\leq H(f_0,f)< 2j\varepsilon_n\bigr\},\alpha\bigr)\Bigr)\leq e^{c_1j^2n\varepsilon_n^2}\,  \Pi(B_{\varepsilon_n^2})^{\alpha}$$
\hskip1.5cm for all sufficiently large positive  integers  $j$  and  $n.$
\smallskip
\noindent Then for each $r_n\to\infty$ we have that $\Pi_n\bigl(A_{r_n\,\varepsilon_n}\bigr)\to 0$ in probability as $n\to\infty$.
\rm
\bigskip
\noindent{\sl Remark.}  Here we adopt the convention that if the denominator of a quotient equals zero then the numerator must also be  zero. Hence, Theorem 3 is still true even when $\Pi(B_{\varepsilon_n^2})=0$ for some $n$.
\bigskip
\noindent {\bf Corollary 4.}  \it Let $\{\varepsilon_n\}_{n=1}^\infty$  be a positive sequence such that  $n\,\varepsilon_n^2$ are uniformly bounded away from zero. Suppose that there exist constants  $c_1$, $c_2$ and a sequence $\{{\cal G}_n\}_{n=1}^\infty$ of subsets on $\F$ such that
$$\leqno \qquad (1)\quad \log N \bigl({j\varepsilon_n\over 3},\bigl\{f\in {\cal G}_n:\,j\varepsilon_n\leq H(f_0,f)< 2j\varepsilon_n\bigr\}\bigr)\leq c_1 n\varepsilon_n^2  $$
\hskip1.5cm for all integer $j$ and $n$ large enough,
$$\leqno \qquad (2)\quad {e^{2n\varepsilon_n^2}\, \Pi(A_{\varepsilon_n}\setminus {\cal G}_n)\over \Pi(B_{\varepsilon_n^2})}\longrightarrow 0 \qquad{\rm as}\quad n\to\infty,$$
$$\leqno \qquad (3)\quad {\Pi\bigl(f\in {\cal G}_n:\,j\varepsilon_n\leq H(f_0,f)< 2j\varepsilon_n\bigr)\over \Pi(B_{\varepsilon_n^2})}\leq e^{c_2j^2n\varepsilon_n^2} \quad {\rm for\ all\ integer\ } j \ {\rm  and }\ n \ {\rm large\ enough}.$$

\noindent Then for each $r_n\to\infty$ we have that $\Pi_n\bigl(A_{r_n\,\varepsilon_n}\bigr)\to 0$ in probability as $n\to\infty$.
\rm
\bigskip
Corollary 4 is a slightly stronger version of Theorem 2.4 of Ghosal et al.(2000). A notable improvement in Corollary 4 is that we have no  restriction on the constant $c_2$, whereas their constant $c_2$ equals half of some universal testing constant.
\bigskip
\noindent{\it Proof of Corollary 4.} We only need to check condition (2) of Theorem 3. It follows from Lemma 1 of Xing et al.(2008) and conditions (1) and (3) that
$$J\bigl({j\varepsilon_n\over 3},\bigl\{f\in {\cal G}_n:\,j\varepsilon_n\leq H(f_0,f)< 2j\varepsilon_n\bigr\},\alpha\bigr)$$$$\leq \alpha\log \Pi \bigl(f\in {\cal G}_n:\,j\varepsilon_n\leq H(f_0,f)< 2j\varepsilon_n\bigr)$$$$+(1-\alpha)\log N\bigl({j\varepsilon_n\over 3},\bigl\{f\in {\cal G}_n:\,j\varepsilon_n\leq H(f_0,f)< 2j\varepsilon_n\bigr\}\bigr)$$
$$\leq \alpha\log\bigl(e^{c_2j^2n\varepsilon_n^2}\Pi(B_{\varepsilon_n^2})\bigr)+(1-\alpha)c_1 n\varepsilon_n^2$$$$=\bigl(\alpha c_2+{(1-\alpha)c_1\over j^2}\bigr)j^2n\varepsilon_n^2+\alpha\log\Pi(B_{\varepsilon_n^2}).$$
Taking a small $\alpha$ in $(0,1)$ and then letting $j$ be large enough, we have that $\alpha c_2+{(1-\alpha)c_1\over j^2}< {1-\alpha\over 18}$ and hence condition (2) of Theorem 3 is fulfilled. The proof of Corollary 4 is complete.
\bigskip
We will conclude this section by presenting an analogue of Theorem 3, which gives an almost sure assertion under stronger conditions.
\bigskip
\noindent {\bf Theorem 4.}  \it Let $\{\varepsilon_n\}_{n=1}^\infty$  be a positive sequence such that there exists a constant $c_0>0$ such that $n\,\varepsilon_n^2\geq c_0\,\log n$ for all large $n$  . Suppose that there exist constants  $0<\alpha<1$, $c_1< {1-\alpha\over 18}$, $c_2>{1\over c_0}$ and a sequence $\{{\cal G}_n\}_{n=1}^\infty$ of subsets on $\F$ such that
$$\leqno \qquad (1)\quad \sum\limits_{n=1}^\infty{e^{n\,\varepsilon_n^2\,(3+2c_2)}\, \Pi(A_{\varepsilon_n}\setminus {\cal G}_n)\over \Pi(W_{\varepsilon_n})}<\infty,$$
$$\leqno \qquad (2)\quad \exp\Bigl( J\bigl({j\varepsilon_n\over 3},\bigl\{f\in {\cal G}_n:\,j\varepsilon_n\leq H(f_0,f)< 2j\varepsilon_n\bigr\},\alpha\bigr)\Bigr)\leq e^{c_1j^2n\varepsilon_n^2}\,  \Pi(W_{\varepsilon_n})^{\alpha}$$
\hskip1.5cm for all sufficiently large positive  integers  $j$  and  $n.$
\smallskip
\noindent Then for each $r$ large enough we have that $\Pi_n\bigl(A_{r\,\varepsilon_n}\bigr)\to 0$ almost surely as $n\to\infty$.
\rm
\bigskip
Completely following the proof of Corollary 4, we have the following consequence of Theorem 4.
\bigskip
\noindent {\bf Corollary 5.}  \it Let $\{\varepsilon_n\}_{n=1}^\infty$  be a positive sequence such that there exists a constant $c_0>0$ such that $n\,\varepsilon_n^2\geq c_0\,\log n$ for all large $n$  . Suppose that there exist constants  $c_1$, $c_2>{1\over c_0}$, $c_3$ and a sequence $\{{\cal G}_n\}_{n=1}^\infty$ of subsets on $\F$ such that
$$\leqno \qquad (1)\quad \log N \bigl({j\varepsilon_n\over 3},\bigl\{f\in {\cal G}_n:\,j\varepsilon_n\leq H(f_0,f)< 2j\varepsilon_n\bigr\}\bigr)\leq c_1 n\varepsilon_n^2  $$
\hskip1.5cm for all integer $j$ and $n$ large enough,
$$\leqno \qquad (2)\quad \sum\limits_{n=1}^\infty{e^{n\,\varepsilon_n^2\,(3+2c_2)}\, \Pi(A_{\varepsilon_n}\setminus {\cal G}_n)\over \Pi(W_{\varepsilon_n})}<\infty,$$
$$\leqno \qquad (3)\quad {\Pi\bigl(f\in {\cal G}_n:\,j\varepsilon_n\leq H(f_0,f)< 2j\varepsilon_n\bigr)\over \Pi(W_{\varepsilon_n})}\leq e^{c_3j^2n\varepsilon_n^2} \quad {\rm for\ all\ integer\ } j \ {\rm  and }\ n \ {\rm large\ enough}.$$

\noindent Then for each $r$ large enough we have that $\Pi_n\bigl(A_{r\,\varepsilon_n}\bigr)\to 0$ almost surely as $n\to\infty$.
\rm
\bigskip
{\bf 3. Illustrations.}\quad In this section we apply our theorems to Bernstein polynomial priors, priors based on uniform distribution, log spline models and finite-dimensional models. This leads some improvements on known results for these models.
\bigskip
\noindent{\it 3.1.  Bernstein polynomial prior.}\quad A Bernstein polynomial prior is a probability measure on the space of continuous probability distribution functions on $[0,1]$.  Petrone (1999) introduced the Bernstein polynomial prior $\Pi$ by putting a prior distribution on the class of Bernstein densities in $[0,1]$ in the following way:
$$b(x;k,F)=\sum\limits_{j=1}^k\bigl(F(j/k)-F((j-1)/k)\bigr)\,\beta(x;j,k-j+1),$$
where $\beta(x;a,b)$ stands for the beta density
$\beta(x;a,b)={\Gamma(a+b)\over \Gamma(a)\Gamma(b)}\,x^{a-1}(1-x)^{b-1},$ $k$ has a distribution $\rho(\cdot)$, $F$ is a random distribution independent of $\rho$. In other words, if $B_j$ stands for the set of all Bernstein densities of order $j$, then $\Pi(\cdot)=\sum_{j=1}^\infty \rho(j)\,\Pi_{B_j}(\cdot)$, where the probability measure $\Pi_{B_j}$ is the normalized restriction of $\Pi$ on $B_j$. We refer to Petrone and Wasserman (2002) for a detailed description of the Bernstein polynomial prior, in which consistency of the posterior distribution for the Bernstein polynomial prior is discussed. Rates of convergence have been established under suitable tail conditions on $\rho$ by Ghosal (2001) and Walker et al.(2007), where the convergence is understood as convergence in $F_0^\infty$-probability.  Ghosal (2001, Theorem 2.3) proved that the prior concentration rate is $(\log n)^{1/3}/n^{1/3}$  and the entropy rate is $(\log n)^{5/6}/n^{1/3}$ under the tail assumption $\rho(j)\approx e^{-cj}$ for all $j$, which yields the convergence rate $(\log n)^{5/6}/n^{1/3}$. Walker et al.(2007) obtained the entropy rate $(\log n)^{1/3}/n^{1/3}$ under the lighter tail condition $\rho(j)\leq e^{-4j\log j}=(1/j^j)^4 $ for all $j$, which leads the convergence rate $(\log n)^{1/3}/n^{1/3}$. In the following we establish the entropy rate $1/n^\gamma$ under the tail condition $\rho(j)\leq (1/j^j)^{c_0} $ for all $j$, where $c_0$ is any fixed positive constant and $\gamma$ is any fixed constant strictly less than $1/2$. Hence we also get the convergence rate $(\log n)^{1/3}/n^{1/3}$ under the weaker condition $\rho(j)\leq (1/j^j)^{c_0} $. This improves the result of Walker et al.(2007).

From Ghosal (2001) it follows that there exists an absolute constant $c>0$ such that $N(\varepsilon_n,B_j)\leq (c/\varepsilon_n)^j$ for all $j$. Given $\gamma<1/2$, choose ${\cal G}_{nj}=B_j$ and $\varepsilon_n=1/n^\gamma$.   Take $0<\alpha<1$ such that $\gamma\,(2+1/d)<1$ with $d=c_0\alpha/(1-\alpha)$. To verify condition (1) of Corollary 3, by $\rho(j)\leq (1/j^j)^{c_0} $ we obtain that
$$ e^{-n\varepsilon_n^2c_1}\,\sum\limits_{j=1}^\infty N(\varepsilon_n, B_j)^{1-\alpha}\,\Pi(B_j)^\alpha\leq e^{-n\varepsilon_n^2c_1}\,\sum\limits_{j=1}^\infty \Bigl({c\over \varepsilon_n}\Bigr)^{j(1-\alpha)}\,\rho(j)^\alpha $$$$\leq e^{-n\varepsilon_n^2c_1}\,\sum\limits_{j=1}^\infty \Bigl({c\over \varepsilon_n\,j^{c_0\alpha\over 1-\alpha}}\Bigr)^{j(1-\alpha)}$$$$\leq e^{-c_1n^{1-2\gamma}}\,\sum\limits_{1\leq j\leq (2cn^\gamma)^{1/d}} \Bigl({cn^\gamma\over j^d}\Bigr)^{j(1-\alpha)}+\sum\limits_{j\geq (2cn^\gamma)^{1/d}} \Bigl({1\over 2}\Bigr)^{j(1-\alpha)},$$
where the last sum on the right hand side tends to zero as $n\to\infty$. To estimate the first term, note that for $g(t)=({cn^\gamma/ t^d})^t$ we have that $g^\prime(t)=({cn^\gamma/ t^d})^t\,\bigl(\log {cn^\gamma/ t^d}-d\bigr)=0$ is equivalent to $t=c^{1/d}v^{\gamma/d}e^{-1}$. This implies $g(j)\leq e^{d_1n^{\gamma/d}}$ for all $j$, where $d_1$ stands for the constant $dc^{1/d}e^{-1}$. Therefore, by the inequality $x\leq e^x$ for $x\geq 0$ we get
$$e^{-c_1n^{1-2\gamma}}\,\sum\limits_{1\leq j\leq (2cn^\gamma)^{1/d}} \Bigl({cn^\gamma\over j^d}\Bigr)^{j(1-\alpha)}\leq e^{-c_1n^{1-2\gamma}}\,(2cn^\gamma)^{1/d}\,e^{d_1(1-\alpha)n^{\gamma/d}}$$
$$\leq e^{-c_1n^{1-2\gamma}+ (2^{1/d}c^{1/d}+d_1-d_1\alpha)n^{\gamma/d}  }\longrightarrow 0\qquad{\rm as}\quad j\to\infty,$$
since $1-2\gamma>\gamma/d$, and hence condition (1) of Corollary 3 holds. Condition (2) of Corollary 3 is trivially fulfilled and hence by Corollary 3 we obtain that the entropy rate is at least $1/n^\gamma$ for any given $\gamma<1/2$.
\bigskip
\noindent{\it 3.2.  Prior based on uniform distribution.}\quad Ghosal et al.(1997) established posterior consistency for prior distributions based on uniform distributions of finite subsets. Priors based on discrete uniform distributions were further studied in Ghosal et al.(2000), in which they used the bracketing entropy as a tool to compute the convergence rate of posterior distributions.
As an application of Theorem 1, we now give a slight extension.  Given  $\varepsilon_n>0$, assume that there exist density functions $f_1,\dots,f_{N_n}$ such that all sets $\{f:\,H_*(f,f_j)\leq 3^{-1/2}\varepsilon_n\}$ form a covering of $\F$. Denote by $\mu_n$ the uniform discrete probability measure on the set $\{f_1,\dots,f_{N_n}\}$. Define a prior distribution $\Pi=\sum_{j=1}^\infty a_j\,\mu_j$ for a given sequence $a_j$ with $a_j>0$ and $\sum_{j=1}^\infty a_j=1$.
\bigskip
\noindent {\bf Theorem 5.}  \it If\quad $\log N_n+\log n+\log {1\over a_n}=O(n\,\varepsilon_n^2)$ as $n\to\infty$, then the posterior distributions $\Pi_n$ converge almost surely at least at the rate $\varepsilon_n$, that is, $\Pi_n\bigl(f:\,H(f_0,f)\geq r\varepsilon_n\bigr)\longrightarrow 0$ as $n\to\infty$ almost surely for any given sufficiently large $r$.
\rm
\bigskip
\noindent{\it Proof.} From $\log n=O(n\,\varepsilon_n^2)$ it follows that $\sum\limits_{n=1}^\infty e^{-n\,\varepsilon_n^2\,c}<\infty$ for all large $c>0$. For any $f$ with $H_*(f,f_j)\leq 3^{-1/2}\varepsilon_n$ we have $H(f,f_j)\leq 3^{1/2}H_*(f,f_j)\leq \varepsilon_n$. Hence we obtain that
$\{f:\,H_*(f,f_j)\leq 3^{-1/2}\varepsilon_n\}\subset
\{f:\,H(f,f_j)<\varepsilon_n\}$ and then $J(\varepsilon_n,\F,0)\leq \log N_n=O(n\,\varepsilon_n^2)$, which implies condition (1) of Theorem 1.  Condition (2) is trivially fulfilled. Condition (3) follows from the fact that $\Pi(f:\,H_*(f_0,f)\leq \varepsilon_n)\geq\Pi(f:\,H_*(f_0,f)\leq 3^{-1/2}\varepsilon_n)\geq {a_n/ N_n}\geq e^{-n\varepsilon_n^2c_3}$ for some $c_3>0$, since $\{f:\,H_*(f_0,f)\leq 3^{-1/2}\varepsilon_n\}$ contains at least some density function of the set $\{f_1,\dots,f_{N_n}\}$. The proof of Theorem 5 is complete.
\bigskip
It seems to be unusual to find a covering of the density space with covering sets of type $\{f:\,H_*(f,f_j)\leq c\,\varepsilon_n\}$. The most widely used norm for continuous functions should be the supremum norm. In fact, one can easily construct a new covering $\cup_{j=1}^{N_n}\{f:\,H_*(f,f_j)\leq c\,\varepsilon_n\}$ of $\F$ in terms of a given covering $\cup_{j=1}^{N_n}\{f:\,||\sqrt{f}-\sqrt{g_j}||_\infty\leq \varepsilon_n\}$ with nonnegative bounded functions $g_j$ (not necessarily density functions), as shown in the following: Take $f_j(x)=(\sqrt{g_j(x)}+\varepsilon_n)^2\big/\int_{\X}(\sqrt{g_j(x)}+\varepsilon_n)^2\,\mu(dx)$, where we assume that $\mu$ is a probability measure on $\X$ and $\varepsilon_n\leq 1$ for all $n$.  Then for each $f$ with $||\sqrt{f}-\sqrt{g_j}||_\infty\leq \varepsilon_n$, that is, $\sqrt{g_j(x)}-\varepsilon_n\leq \sqrt{f(x)}\leq \sqrt{g_j(x)}+ \varepsilon_n$ on $\X$, we have
$${f(x)\over f_j(x)}={f(x)\over (\sqrt{g_j(x)}+ \varepsilon_n)^2}\,\int_{\X} \bigl(\sqrt{g_j(x)}+ \varepsilon_n\bigr)^2\,\mu(dx) \leq \int_{\X} \bigl(\sqrt{f^*_j(x)}+2\, \varepsilon_n\bigr)^2\,\mu(dx)$$
$$=1+4\,\varepsilon_n^2+ 4\,\varepsilon_n\,\int_{\X}\sqrt{f^*_j(x)} \,\mu(dx)\leq (1+2\,\varepsilon_n)^2,$$
where $f^*_j$ is some density function in $\{f:\,||\sqrt{f}-\sqrt{g_j}||_\infty\leq \varepsilon_n\}$ and the last inequality follows from $||\sqrt{f_*}||_1\leq ||\sqrt{f_*}||_2=1$. This implies that
$$H_*(f,f_j)\leq \Bigl(\,{2\over 3}\,(1+2\,\varepsilon_n)+{1\over 3}\,\Bigr)^{1\over2}\,H(f,f_j)\leq 2\,H(f,f_j)$$$$\leq  2\,H\big(f,(\sqrt{g_j}+\varepsilon_n)^2\big)+2\,H\big((\sqrt{g_j}+\varepsilon_n)^2,f_j\big)$$$$\leq 4\,\varepsilon_n+
2\,\biggl(\Bigl(\int_{\X} \bigl(\sqrt{g_j(x)}+ \varepsilon_n\bigr)^2\,\mu(dx)\Bigr)^{1\over 2}-1\biggr)\leq 4\,\varepsilon_n+
2\,\bigl((1+2\,\varepsilon_n)-1\bigr)=8\,\varepsilon_n.$$
Therefore, we have
 $$\bigcup\limits_{j=1}^{N_n}\big\{f:\,H_*(f,f_j)\leq 8\,\varepsilon_n\big\}\supset \bigcup_{j=1}^{N_n}\big\{f:\,||f-g_j||_\infty\leq \varepsilon_n\big\}\supset \F.$$
Observe that the numbers of covering subsets in both type coverings are equal.

For models with $H_*(f,g)$ controlled by a constant multiple of the Hellinger metric $H(f,g)$ such as the exponential family and  a model with uniformly bounded supremum norm $||f/g||_\infty$ for all density functions $f$ and $g$, it is not necessary to assume that the probability measures $\mu_n$ constructed above concentrate on a finite number of points. Here we give an extension of Theorem 5. Let $c_0\geq1$ and $\F_{c_0}$ be a subfamily of $\F$ such that $H_*(f,g)\leq c_0\, H(f,g)$ for all $f,g\in \F_{c_0}$. Given  $\varepsilon_n>0$, let $\{P_1,\dots,P_{K_n}\}$ be a partition of $\F_{c_0}$ such that for each $P_i$ there exists $f_i$ in $\F_{c_0}$ with $P_i\subset \{f:\, H(f_i,f)\leq \varepsilon_n/2\,c_0\}$. Take any probability measure $\bar\mu_n$ on $\F_{c_0}$ with $\bar\mu_n(P_i)=1/K_n$ for $i=1,2,\dots,K_n$. Define then a prior distribution $\bar\Pi=\sum_{j=1}^\infty a_j\,\bar \mu_j$ for a given sequence $a_j$ with $a_j>0$ and $\sum_{j=1}^\infty a_j=1$. Now we have
\bigskip
\noindent {\bf Theorem 6.}  \it Let $f_0\in \F_{c_0}$. If\quad $\log K_n+\log n+\log {1\over a_n}=O(n\,\varepsilon_n^2)$ as $n\to\infty$, then the posterior distributions $\bar\Pi_n$ converge at least at the rate $\varepsilon_n$ almost surely.
\rm
\bigskip
\noindent{\it Proof.} By the proof of Theorem 5, we only need to verify condition (3) of Theorem 1. Take $f_{i_0}\in \F_{c_0}$ such that $H(f_0,f_{i_0})\leq \varepsilon_n/2\,c_0$.
Then, for all $f\in \F_{c_0}$ with $H(f_{i_0},f)\leq \varepsilon_n/2\,c_0$, we have that $H_*(f_0,f)\leq c_0\,H(f_0,f)\leq c_0\,H(f_0,f_{i_0})+c_0\,H(f_{i_0},f)\leq \varepsilon_n. $ Hence we get that
$\Pi(f:\,H_*(f_0,f)\leq \varepsilon_n)\geq\Pi(P_{i_0})\geq {a_n/ K_n}\geq e^{-n\varepsilon_n^2c_3}$ for some $c_3>0$, and the proof of Theorem 6 is complete.
\bigskip
Observe that, given a covering $\{O_1,O_2,\dots,O_{K_n}\}$ of $\F_{c_0}$, one can easily construct a partition $\{P_1,P_2,\dots,P_{K_n}\}$ of $\F_{c_0}$ in the following way: $P_1=O_1\cap \F_{c_0}$ and $P_i=(O_i-\cup_{l=1}^{i-1}P_l)\cap \F_{c_0}$ for $i=2,3,\dots,K_n$.
\bigskip
\noindent {\it Example} ( Exponential families ).  We consider the exponential family of all density functions of the form $e^{h(x)}$, where the function $h(x)$ belongs to a fixed bounded subset in the Sobolev space $C^p[0,1]$ with $p>0$. A subclass of this family has been recently studied by Scricciolo (2006). Following a result of Kolmogorov and Tihomirov (1959), we know that the $\varepsilon$-entropy of this family with respect to the norm $||\cdot||_\infty$ equals $O(\varepsilon^{-1/p})$. Thus, using the above argument we get that $\log K_n=O(n\,\varepsilon_n^2)$ for $\varepsilon_n=n^{-p/(2p+1)}$ and hence by Theorem 6 the posterior distributions constructed above converge at the rate $\varepsilon_n=n^{-p/(2p+1)}$, which is known to be the optimal rate of convergence in the minimax sense under the Hellinger loss.
\bigskip
\noindent{\it 3.3.  Log spline models.}\quad Log spline models for density estimation have been studied, among others, by Stone (1990) and Ghosal et al.(2000). Let $\big[ (k-1)/K_n,k/K_n\big)$ with $k=1,2,\dots,K_n$ be a partition of the interval $[0,1)$. The space of splines of order $q$ relative to this partition is the set of all functions $f:[0,1)\to \R$ such that $f$ is $q-2$ times continuously differentiable on $[0,1)$ and the restriction of $f$ on each $\big[ (k-1)/K_n,k/K_n\big)$ is a polynomial of degree strictly less then $q$. Let $J_n=q+K_n-1$. This space of splines is a $J_n$-dimensional vector space with a B-spline basis $B_1(x),B_2(x),\dots,B_{J_n}(x)$, see Ghosal et al.(2000) for the details of such a basis. Let $\F$ be the set of all density functions in $C^\alpha[0,1]$. Assume that the true density function $f_0(x)$ is bounded away from zero and infinity. We consider the $J_n$-dimensional exponential  subfamily of $C^\alpha[0,1]$ of the form
$$f_\theta(x)=\exp\Big(\ \sum\limits_{j=1}^{J_n}\theta_jB_j(x)-c(\theta)\ \Big),$$
where $\theta=(\theta_1,\theta_2,\dots,\theta_{J_n})\in \Theta_0=\{(\theta_1,\theta_2,\dots,\theta_{J_n})\in \R^{J_n}: \sum_{j=1}^{J_n}\theta_j=0\}$ and the constant $c(\theta)$ is chosen such that $f_\theta(x)$ is a density function in $[0,1]$. Each prior on $\Theta_0$ induces naturally a prior on $\F$. Let $||\theta||_\infty=\max_j|\theta_j|$ be the infinity norm on $\Theta_0$. Assume that $a_1\,n^{1/(2\alpha+1)}\leq K_n\leq a_2\,n^{1/(2\alpha+1)}$ for two fixed positive constants $a_1$ and $a_2$. Assume that the prior $\Pi$ for $\Theta_0$ is supported on $[-M,M]^{J_n}$ for some $M\geq 1$ and has a density function with respect to the Lebesgue measure on $\Theta_0$, which is bounded below by $a_3^{J_n}$ and above by  $a_4^{J_n}$. Take a constant $d>0$ such that $d\,||\theta||_\infty\leq ||\log f_\theta(x)||_\infty$ for all $\theta\in\Theta_0$. Ghosal et al.(2000, Theorem 4.5) proved that, if $f_0\in C^\alpha[0,1]$ with $q\geq \alpha\geq 1/2$ and $||\log f_0(x)||_\infty\leq d\,M/2,$ the posteriors $\Pi_n$ converge in probability at the rate $n^{-\alpha/(2\alpha+1)}$. Using Corollary 5 we now get that under the same assumptions as in Ghosal et al.(2000, Theorem 4.5), the posteriors $\Pi_n$ are in fact convergent almost surely at the rate $\varepsilon_n=n^{-\alpha/(2\alpha+1)}$. To see this, take ${\cal G}_n=\F$. Clearly, $n\,\varepsilon_n^2\geq \log n$ for all large $n$. Condition (1) of Corollary 5 has been verified by Ghosal et al.(2000) and condition (2) is trivially fulfilled. Condition (3) follows also from the proof of Theorem 4.5 in  Ghosal et al.(2000), since the inequality $H_*(f_0,f_\theta)\leq H(f_0,f_\theta)\,\big|\big|{f_0/ f_\theta}\big|\big|_\infty^{1/2}$ holds for all $\theta\in \Theta_0$.
\bigskip
\noindent{\it 3.4.  Finite-dimensional models.}\quad Let $\beta>0$ and let $\Theta$ be a bounded subset in $\R^d$ with the Euclidean norm $||\cdot||$. Denote by $\F$ the family of all density functions $f_\theta$ with the parameter $\theta$ in $\Theta$ satisfying
$$a_1\,||\theta_1-\theta_2||^\beta\leq H(f_{\theta_1},f_{\theta_2})\leq \sqrt{3}\,H_*(f_{\theta_1},f_{\theta_2})\leq a_2\,||\theta_1-\theta_2||^\beta$$
for all $\theta_1,\, \theta_2\in \Theta$, where $a_1$ and $a_2$ are two fixed positive constants. Assume that the true value $\theta_0$ is in $ \Theta$ and that the density function of the prior distribution $\Pi$ with respect to the Lebesgue measure on $\Theta$ is uniformly bounded away from zero and infinity. Under slightly weaker conditions, Ghosal et al.(2000) proved that the posterior distributions $\Pi_n$ converge in probability at the rate $1/\sqrt n$.
Now we give an almost sure assertion for this model.
\bigskip
\noindent {\bf Theorem 7.}  \it Under the above assumptions,  the posterior distributions $\Pi_n$ converge almost surely at least at the rate $\sqrt{\log n}/\sqrt{n}$.
\rm
\bigskip
\noindent{\it Proof.} We shall apply Corollary 5 for ${\cal G}_n=\F$. Clearly, $n\,\varepsilon_n^2=\,\log n$ for $\varepsilon_n=\sqrt{\log n}/\sqrt{n}$. Condition (1) has been verified in the proof of Theorem 5.1 of Ghosal et al.(2000). Condition (2) is trivially fulfilled. Using $H_*(f_{\theta_0},f_{\theta})\leq 3^{-1/2}\,a_2\,||\theta_0-\theta||^\beta$, we have that $\Pi(W_{\varepsilon_n})\geq \Pi\big(\theta:\,||\theta-\theta_0||\leq (\sqrt{3}\,\varepsilon_n/a_2)^{1/\beta}\big)$. Hence, the verification of condition (3) follows from the same lines as the proof of Ghosal et al.(2000, Theorem 5) and then by Corollary 5 we conclude the proof of Theorem 7.
\bigskip\bigskip
{\bf 4. Lemmas and Proofs.}\quad In this section we give proofs of our lemmas and theorems. For simplicity of notations, we assume throughout this section that $\varepsilon_n=\bar\varepsilon_n=\tilde\varepsilon_n$.
\bigskip
\noindent{\it Proof of Lemma 1.} It is no restriction to assume that $\Pi\bigl( W_{\varepsilon})>0$. Using Jensen's inequality for the convex function $x^{-1/2}$ for $x>0$ and Chebyshev's inequality, we obtain that
$$F^\infty_0\Bigl(\ \int_{\F} R_n(f)\, \Pi(df)\leq e^{-n\,\varepsilon^2\,(3+2c)}\ \Pi\bigl( W_{\varepsilon})\, \Bigr)$$$$\leq F^\infty_0\Bigl(\ \int_{W_{\varepsilon}} R_n(f)\, \Pi(df)\leq e^{-n\,\varepsilon^2\,(3+2c)}\ \Pi\bigl( W_{\varepsilon})\, \Bigr)$$
$$=F^\infty_0\biggl(\ e^{n\,\varepsilon^2\,({3\over 2}+c)}\leq  \Bigl(\,{1\over \Pi(W_{\varepsilon}\,)}\,\int_{W_{\varepsilon}} R_n(f)\, \Pi(df)\Bigr)^{-{1\over 2}}\ \biggr)$$
$$\leq F^\infty_0\biggl(\ e^{n\,\varepsilon^2\,({3\over 2}+c)}\leq {1\over \Pi(W_{\varepsilon}\,)}\,\int_{W_{\varepsilon}} R_n(f)^{-{1\over 2}}\, \Pi(df)\ \biggr)$$
$$\leq e^{-n\,\varepsilon^2\,({3\over 2}+c)}\, {1\over \Pi(W_{\varepsilon}\,)}\,E\,\int_{W_{\varepsilon}} R_n(f)^{-{1\over 2}}\, \Pi(df)$$$$=e^{-n\,\varepsilon^2\,({3\over 2}+c)}\,{1\over \Pi(W_{\varepsilon}\,)}\,\int_{W_{\varepsilon}} \biggl(\, E\, \sqrt{f_0(X_1)\over f(X_1)}\ \biggr)^n\, \Pi(df).$$
On the other hand, we have
$$E\, \sqrt{f_0(X_1)\over f(X_1)}=1+E\,{\sqrt{f_0(X_1)}-\sqrt{f(X_1)}\over \sqrt{f(X_1)}} $$
$$=1+\int_{\X} {\sqrt{f_0(x)}-\sqrt{f(x)}\over \sqrt{f(x)}}\,\Bigl(\, \sqrt{f_0(x)}-\sqrt{f(x)}+\sqrt{f(x)}\,\Bigr) \,\sqrt{f_0(x)}\, \mu(dx)$$
$$=1+\int_{\X} \Bigl(\, \sqrt{f_0(x)}-\sqrt{f(x)}\,\Bigr)^2 \,{\sqrt{f_0(x)}\over \sqrt{f(x)}} \, \mu(dx)+\int_{\X} \bigl(\, f_0(x)-\sqrt{f(x)\,f_0(x)}\,\bigr)\,  \mu(dx)$$
$$= 1+ \int_{\X} \Bigl(\, \sqrt{f_0(x)}-\sqrt{f(x)}\,\Bigr)^2 \,{\sqrt{f_0(x)}\over \sqrt{f(x)}} \, \mu(dx)+{1\over 2}\,\int_{\X} \bigl(\, \sqrt{f_0(x)}-\sqrt{f(x)}\,\bigr)^2\,  \mu(dx)$$
$$= 1+ {3\over 2}\,H_*(f_0,f)^2\leq e^{{3\over 2}\,H_*(f_0,f)^2 }\leq e^{{3\over 2}\,\varepsilon^2}, $$
where the last inequality holds when $f\in W_{\varepsilon} $. Hence we get
$$F^\infty_0\Bigl(\ \int_{\F} R_n(f)\, \Pi(df)\leq e^{-n\,\varepsilon^2\,(3+2c)}\ \Pi\bigl( W_{\varepsilon})\, \Bigr)$$$$\leq e^{-n\,\varepsilon^2\,({3\over 2}+c)}\, {1\over \Pi(W_{\varepsilon}\,)}\,\int_{W_{\varepsilon}} e^{{3\over 2}\,n\,\varepsilon^2}  \, \Pi(df)= e^{-n\,\varepsilon^2\,c}.$$
The proof of Lemma 1 is complete.
\bigskip
In the proof of Theorem 1 we use the following Lemma, which is similar to Lemma 5 of Barron et al.(1999).
\bigskip
\noindent {\bf Lemma 2.}  \it Let $c_2> 0$ and $c_3\geq 0$. Let $\{\varepsilon_n\}_{n=1}^\infty$  be a positive sequence such that  $\Pi(W_{\varepsilon_n})\geq e^{-n\,\varepsilon_n^2\,c_3}$ for all $n $ and $\sum\limits_{n=1}^\infty e^{-n\,\varepsilon_n^2\,c_2}<\infty$. If a sequence $\{D_n\}_{n=1}^\infty$ of subsets in $\F$ satisfies
$\sum\limits_{n=1}^\infty e^{n\,\varepsilon_n^2\,(3+3c_2+c_3)}\,\Pi(D_n)<\infty$, then $\Pi_n(D_n)\to 0$ almost surely as $n\to\infty$.
\rm
\bigskip
\noindent{\it Proof.} From Chebyshev's inequality and Fubini's theorem it turns out that
$$F_0^\infty\Bigl\{\ \int_{D_n}R_n(f)\,\Pi(df)\geq e^{-n\,\varepsilon_n^2\,(3+3c_2+c_3)}\ \Bigr\}\leq e^{n\,\varepsilon_n^2\,(3+3c_2+c_3)}\ E\int_{D_n}R_n(f)\,\Pi(df)$$
$$=e^{n\,\varepsilon_n^2\,(3+3c_2+c_3)}\, \int_{D_n}E\,R_n(f)\,\Pi(df)=e^{n\,\varepsilon_n^2\,(3+3c_2+c_3)}\,\Pi(D_n)$$
for all $n$. Hence by the first Borel-Cantelli Lemma we get that
$$\int_{D_n}R_n(f)\,\Pi(df)\leq e^{-n\,\varepsilon_n^2\,(3+3c_2+c_3)}$$
almost surely for all $n$ large enough. On the other hand, Lemma 1 and the first Borel-Cantelli Lemma yield that
$$\int_{\F} R_n(f)\, \Pi(df)\geq \Pi(W_{\varepsilon_n})\,e^{-n\varepsilon_n^2(3+2c_2)}\geq  e^{-n\,\varepsilon_n^2\,(3+2c_2+c_3)}$$
almost surely for all $n$. Therefore, we obtain that with probability one,
 $$\Pi_n(D_n)={\int_{D_n} R_n(f)\, \Pi(df)\over \int_{\F} R_n(f)\, \Pi(df)}\leq e^{-n\,\varepsilon_n^2\,c_2},$$
which tends to zero as $n\to\infty$ and the proof of Lemma 2 is complete.
\bigskip
\noindent{\it Proof of Theorem 1.} It is clear that if condition (1) holds for some $\alpha=\alpha_0$ then it also holds for any $\alpha\geq\alpha_0$. So we may assume that $0<\alpha<1$.  Given $r>2+\sqrt{2(3\alpha +2\alpha c_2+\alpha c_3+c_1)\over 1-\alpha}$, we have
$$\Pi_n(A_{r\varepsilon_n})\leq \Pi_n\bigl({\cal G}_n\cap A_{r\varepsilon_n}\bigr)+\Pi_n\bigl(A_{\varepsilon_n}\setminus {\cal G}_n\bigr).$$
It then follows from Lemma 2 that $\Pi_n\bigl(A_{\varepsilon_n}\setminus {\cal G}_n\bigr)\to 0$ almost surely as $n\to\infty$. So it suffices to prove that $\Pi_n\bigl({\cal G}_n\cap A_{r\varepsilon_n}\bigr)\to 0$ almost surely as $n\to\infty$.
By the definition of $J(\varepsilon_n,{\cal G}_n,\alpha)$, for each fixed $n$ there exist functions $f_1,f_2,\dots,f_N$ in $L_\mu$ such that
${\cal G}_n\cap A_{r\varepsilon_n}\subset\bigcup_{j=1}^NB_j$, where $B_j= {\cal G}_n\cap A_{r\varepsilon_n}\cap\{f:\,H(f_j,\,f)<\varepsilon_n\}$ and $\sum_{j=1}^N \Pi(B_j)^\alpha\leq 2\,e^{J(\varepsilon_n,{\cal G}_n,\alpha)}$.
It is no restriction to assume that all the sets $B_j$ are disjoint and nonempty. Taking a $f_j^\star\in B_j$ we get that
$H(f_j,f_0)\geq H(f_j^\star,f_0)-H(f_j^\star,f_j)\geq (r-1)\,\varepsilon_n$.
Now for each $B_j$ we have
$$\int_{B_j}R_n(f)\, \Pi(df)=\Pi({B_j})\,\prod\limits_{k=0}^{n-1}\,{\int_{B_j}R_{k+1}(f)\, \Pi(df)\over \int_{B_j}R_k(f)\, \Pi(df)}=\Pi({B_j})\,\prod\limits_{k=0}^{n-1}\,{f_{k{B_j}}(X_{k+1})\over f_0(X_{k+1})},$$
where
$f_{k{B_j}}(x)={\int_{B_j}f(x)\,R_k(f)\, \Pi(df)\big/ \int_{B_j}R_k(f)\, \Pi(df)}$ and $R_0(f)=1$.
The function $f_{k{B_j}}$ was introduced by Walker (2004) and can be considered  as the predictive density of $f$ with a normalized posterior distribution, restricted on the set ${B_j}$.
Clearly, Jensen's inequality yields that $H(f_{kB_j},f_j)^2\leq \varepsilon_n^2$ for each $k$. Hence
$H(f_{kB_j},f_0)\geq H(f_j,f_0)-H(f_j,f_{kB_j})\geq (r-2)\,\varepsilon_n>0$.
Since
$\sum\limits_{n=1}^\infty e^{-n\,\varepsilon_n^2\,c_2}<\infty$, it turns out from Lemma 1 and the first Borel-Cantelli Lemma that
$\int_{\F} R_n(f)\, \Pi(df)\geq e^{-n\,\varepsilon_n^2\,(3+2c_2)}\ \Pi\bigl( W_{\varepsilon_n})$ almost surely for all $n$ large enough.
Hence, by condition (3) we obtain that
$$\Pi_n\bigl({\cal G}_n\cap A_{r\varepsilon_n}\bigr)\leq \bigl(\Pi_n({\cal G}_n\cap A_{r\varepsilon_n})\bigr)^\alpha\leq \bigl(\ \sum\limits_{j=1}^N\Pi_n(B_j)\ \bigr)^\alpha$$
$$ \leq \sum\limits_{j=1}^N\Pi_n(B_j)^\alpha={{\sum\limits_{j=1}^N \Pi(B_j)^\alpha\,\prod\limits_{k=0}^{n-1}\,{f_{kB_j}(X_{k+1})^\alpha\over f_0(X_{k+1})^\alpha}}\over \Bigl(\int_{\F} R_n(f)\, \Pi(df) \Bigr)^\alpha }\leq   {{\sum\limits_{j=1}^N \Pi(B_j)^\alpha\,\prod\limits_{k=0}^{n-1}\,{f_{kB_j}(X_{k+1})^\alpha\over f_0(X_{k+1})^\alpha}}\over \Bigl(\Pi(W_{\varepsilon_n})\,e^{-n\varepsilon_n^2(3+2c_2)} \Bigr)^\alpha }$$
$$\leq  e^{n\,\varepsilon_n^2\,(3+2c_2+c_3)\,\alpha}  {\sum\limits_{j=1}^N \Pi(B_j)^\alpha\,\prod\limits_{k=0}^{n-1}\,{f_{kB_j}(X_{k+1})^\alpha\over f_0(X_{k+1})^\alpha}}$$
almost surely for all $n$ large enough.
Since $r>2+\sqrt{2(3\alpha+2\alpha c_2+\alpha c_3+c_1)\over 1-\alpha}$, the inequality $(3+2 c_2+c_3)\,\alpha<{1\over 2}\,(r-2)^2\,(1-\alpha)-c_1$ holds. Take a constant $b$ with $(3+2 c_2+c_3)\,\alpha<b<{1\over 2}\,(r-2)^2\,(1-\alpha)-c_1$. Denote ${\cal F}_k=\sigma\{X_1,X_2,\dots,X_k\}$. Then we have
$$F^\infty_0\Bigl\{\sum\limits_{j=1}^N \Pi(B_j)^\alpha\,\prod\limits_{k=0}^{n-1}\,{f_{kB_j}(X_{k+1})^\alpha\over f_0(X_{k+1})^\alpha}\geq e^{-n\,\varepsilon_n^2\,b}\Bigr\}$$$$\leq e^{n\,\varepsilon_n^2\,b}\,E\,\biggl(\,\sum\limits_{j=1}^N \Pi(B_j)^\alpha\,\prod\limits_{k=0}^{n-1}\,{f_{kB_j}(X_{k+1})^\alpha\over f_0(X_{k+1})^\alpha}\,\biggr)$$
$$=e^{n\,\varepsilon_n^2\,b}\,\sum\limits_{j=1}^N \Pi(B_j)^\alpha\,E\,\biggl(\,\prod\limits_{k=0}^{n-1}\,{f_{kB_j}(X_{k+1})^\alpha\over f_0(X_{k+1})^\alpha}\,\biggr)$$
$$\leq 2\,e^{J(\varepsilon_n,{\cal G}_n,\alpha)+n\,\varepsilon_n^2\,b}\,\max\limits_{1\leq j\leq N}\,E\,\biggl(\,\prod\limits_{k=0}^{n-1}\,{f_{kB_j}(X_{k+1})^\alpha\over f_0(X_{k+1})^\alpha}\,\biggr),$$
where
$$E\,\biggl(\,\prod\limits_{k=0}^{n-1}\,{f_{kB_j}(X_{k+1})^\alpha\over f_0(X_{k+1})^\alpha}\,\biggr)=E\,\Biggl(E\biggl(\ \prod\limits_{k=0}^{n-1}\,{f_{kB_j}(X_{k+1})^\alpha\over f_0(X_{k+1})^\alpha}\ \bigg|\ {\cal F}_{n-1}\ \biggr)\Biggr) $$
$$=E\Biggl(\prod\limits_{k=0}^{n-2}\,{f_{kB_j}(X_{k+1})^\alpha\over f_0(X_{k+1})^\alpha}\ E\biggl(\ {f_{n-1B_j}(X_n)^\alpha\over f_0(X_n)^\alpha}\ \bigg|\ {\cal F}_{n-1}\ \biggr)\Biggr). $$
By the conditional H\"older's inequality we get that with probability one,
$$E\biggl(\ {f_{n-1B_j}(X_n)^\alpha\over f_0(X_n)^\alpha}\ \bigg|\ {\cal F}_{n-1}\ \biggr)=E\biggl(\ {f_{n-1B_j}(X_n)^{\alpha\over 2}\over f_0(X_n)^{\alpha\over 2}}\ {f_{n-1B_j}(X_n)^{\alpha\over 2}\over f_0(X_n)^{\alpha\over 2}}\ \bigg|\ {\cal F}_{n-1}\ \biggr)$$
$$\leq E\biggl(\ {f_{n-1B_j}(X_n)^{{\alpha\over 2}\cdot{2\over 2-\alpha}}\over f_0(X_n)^{{\alpha\over 2}\cdot{2\over 2-\alpha}}}\ \bigg|\ {\cal F}_{n-1}\ \biggr)^{2-\alpha\over 2}\ E\biggl(\ {f_{n-1B_j}(X_n)^{{\alpha\over 2}\cdot{2\over \alpha}}\over f_0(X_n)^{{\alpha\over 2}\cdot{2\over \alpha}}}\ \bigg|\ {\cal F}_{n-1}\ \biggr)^{\alpha\over 2}$$
$$= E\biggl(\ {f_{n-1B_j}(X_n)^{\alpha\over 2-\alpha}\over f_0(X_n)^{\alpha\over 2-\alpha}}\ \bigg|\ {\cal F}_{n-1}\ \biggr)^{2-\alpha\over 2}.$$
Take the integer $m$ with ${\alpha\over 1-\alpha}\leq 2^m<{2\alpha\over 1-\alpha}.$
Repeating the above procedure $m-1$ more times we obtain that with probability one,
$$E\biggl(\ {f_{n-1B_j}(X_n)^\alpha\over f_0(X_n)^\alpha}\ \bigg|\ {\cal F}_{n-1}\ \biggr)\leq E\biggl(\ {f_{n-1B_j}(X_n)^{\alpha\over 2^m(1-\alpha)+\alpha}\over f_0(X_n)^{\alpha\over 2^m(1-\alpha)+\alpha}}\ \bigg|\ {\cal F}_{n-1}\ \biggr)^{2^m(1-\alpha)+\alpha\over 2^m},$$
which by the conditional H\"older's inequality is less than
$$E\biggl(\ {f_{n-1B_j}(X_n)^{1\over 2}\over f_0(X_n)^{1\over 2}}\ \bigg|\ {\cal F}_{n-1}\ \biggr)^{\alpha\over 2^{m-1}}=\Bigl(\ \int\sqrt{f_{n-1B_j}(X_n)\ f_0(X_n)}\ \mu(dX_n) \Bigr)^{\alpha\over 2^{m-1}}$$$$
= \Bigl(1-{H(f_{n-1B_j},f_0)^2\over 2}\Bigr)^{\alpha\over 2^{m-1}}
\leq \Bigl(1-{(r-2)^2\,\varepsilon_n^2\over 2}\Bigr)^{\alpha\over 2^{m-1}}$$$$\leq e^{-2^{-m}\,(r-2)^2\,\alpha\,\varepsilon_n^2}\leq e^{{1\over 2}\,(r-2)^2\,(\alpha-1)\,\varepsilon_n^2}.$$
Hence, with probability one, we have
$$E\biggl(\ \prod\limits_{k=0}^{n-1}\,{f_{kB_j}(X_{k+1})^\alpha\over f_0(X_{k+1})^\alpha}\ \biggr)\leq e^{{1\over 2}\, (r-2)^2\,(\alpha-1)\,\varepsilon_n^2}\ E\Biggl(\ \prod\limits_{k=0}^{n-2}\,{f_{kB_j}(X_{k+1})^\alpha\over f_0(X_{k+1})^\alpha}\ \Biggr). $$
Repeating the same argument $n-1$ times, we obtain that for each $j$,
$$E\biggl(\ \prod\limits_{k=0}^{n-1}\,{f_{kB_j}(X_{k+1})^\alpha\over f_0(X_{k+1})^\alpha}\ \biggr)\leq e^{{1\over 2}\,(r-2)^2\,(\alpha-1)\,n\,\varepsilon_n^2}. $$
Therefore, we have gotten that for all $n,$
$$F^\infty_0\Bigl\{\,\sum\limits_{j=1}^N \Pi(B_j)^\alpha\,\prod\limits_{k=0}^{n-1}\,{f_{kB_j}(X_{k+1})^\alpha\over f_0(X_{k+1})^\alpha}\geq e^{-n\,\varepsilon_n^2\,b}\,\Bigr\}$$$$\leq 2\,e^{J(\varepsilon_n,{\cal G}_n,\alpha)+n\,\varepsilon_n^2\bigl(b+{1\over 2}\,(r-2)^2\,(\alpha-1)\bigr)}\leq 2\, e^{J(\varepsilon_n,{\cal G}_n,\alpha)-n\,\varepsilon_n^2\,c_1}.$$
Thus, together with condition (1), the first Borel-Cantelli Lemma yields that
$$\sum\limits_{j=1}^N \Pi(B_j)^\alpha\,\prod\limits_{k=0}^{n-1}\,{f_{kB_j}(X_{k+1})^\alpha\over f_0(X_{k+1})^\alpha}\leq e^{-n\,\varepsilon_n^2\,b}$$
almost surely for all $n$ large enough. Hence we have
$$\Pi_n\bigl({\cal G}_n\cap A_{r\varepsilon_n}\bigr)\leq 2\,e^{n\,\varepsilon_n^2\,(3\,\alpha+2\,\alpha\,c_2+\alpha\,c_3-b)},  $$
which tends to zero as $n\to\infty$, since $(3+2 c_2+c_3)\,\alpha<b$ and $n\,\varepsilon_n^2\to \infty$ as $n\to\infty$. The proof of Theorem 1 is complete.
\bigskip
To prove Theorem 2, we need a replacement of Lemma 2 under weaker conditions.
\bigskip
\noindent {\bf Lemma 3.}  \it Let $c_2\geq 0$ and let $\{\varepsilon_n\}_{n=1}^\infty$  be a positive sequence such that  $\Pi(B_{\varepsilon_n^2})\geq e^{-n\,\varepsilon_n^2\,c_2}$ for all $n $. If a sequence $\{D_n\}_{n=1}^\infty$ of subsets in $\F$ satisfies
$e^{n\,\varepsilon_n^2\,(2+c_2)}\,\Pi(D_n)\to 0$ as $n\to\infty$, then $\Pi_n(D_n)\to 0$ in probability as $n\to\infty$.
\rm
\bigskip
\noindent{\it Proof.} From Lemma 1 of Shen et al. (2001) or Lemma 8.1 of Ghosal et al. (2000)
it turns out that we have, with probability tending to 1,
 $$\Pi_n(D_n)\leq{\int_{D_n} R_n(f)\, \Pi(df)\over \Pi(B_{\varepsilon_n^2})\,e^{-2n \varepsilon_n^2}}\leq e^{n\,\varepsilon_n^2\,(2+c_2)}\,\int_{D_n} R_n(f)\, \Pi(df).$$
Hence for any given $\delta>0$ we have that
$$F_0^\infty\bigl\{\Pi_n(D_n)\geq \delta\bigr\}\leq F_0^\infty\Bigl\{ e^{n\,\varepsilon_n^2\,(2+c_2)}\,\int_{D_n} R_n(f)\, \Pi(df)  \geq \delta\Bigr\}+{\rm o}(1)$$
$$\leq  {1\over \delta}\,e^{n\,\varepsilon_n^2\,(2+c_2)}\,E\,\int_{D_n} R_n(f)\, \Pi(df)+{\rm o}(1)$$$$={1\over \delta}\,e^{n\,\varepsilon_n^2\,(2+c_2)}\,\Pi(D_n)+{\rm o}(1)\longrightarrow 0\qquad {\rm as}\quad n\to\infty, $$
which concludes the proof of Lemma 3.
\bigskip
\noindent{\it Proof of Theorem 2.} Assume that $0<\alpha<1$. The proof of Theorem 2 follows from the same lines as the proof of Theorem 1.
By Lemma 3 it suffices to prove that $\Pi_n\bigl({\cal G}_n\cap A_{r\,\varepsilon_n}\bigr)\to 0$ in probability as $n\to\infty$. For any given $\delta>0$, following the proof of Theorem 1 we get
$$F_0^\infty\bigl\{\Pi_n\bigl({\cal G}_n\cap A_{r\,\varepsilon_n}\bigr)\geq \delta\bigr\}\leq F_0^\infty\biggl\{
e^{n\,\varepsilon_n^2\,(2+c_2)\,\alpha}  {\sum\limits_{j=1}^N \Pi(B_j)^\alpha\,\prod\limits_{k=0}^{n-1}\,{f_{kB_j}(X_{k+1})^\alpha\over f_0(X_{k+1})^\alpha}}
\geq \delta\biggr\}+{\rm o}(1)$$
$$\leq {1\over \delta}\,e^{n\,\varepsilon_n^2\,(2+c_2)\,\alpha}\, {\sum\limits_{j=1}^N \Pi(B_j)^\alpha\,E\,\prod\limits_{k=0}^{n-1}\,{f_{kB_j}(X_{k+1})^\alpha\over f_0(X_{k+1})^\alpha}}+{\rm o}(1)$$
$$\leq {2\over \delta}\,e^{J(\varepsilon_n,{\cal G}_n,\alpha)+n\,\varepsilon_n^2\,(2+c_2)\,\alpha}
\,\max\limits_{1\leq j\leq N}\,E\,\prod\limits_{k=0}^{n-1}\,{f_{kB_j}(X_{k+1})^\alpha\over f_0(X_{k+1})^\alpha}+{\rm o}(1)$$
$$\leq {2\over \delta}\,e^{J(\varepsilon_n,{\cal G}_n,\alpha)+n\,\varepsilon_n^2\,\bigl((2+c_2)\,\alpha+{1\over 2}(r-2)^2(\alpha-1)\bigr)}+{\rm o}(1)$$$$\leq {2\over \delta}\,e^{J(\varepsilon_n,{\cal G}_n,\alpha)-n\,\varepsilon_n^2\,c_2}+{\rm o}(1)\longrightarrow 0\qquad {\rm as}\quad n\to\infty,$$
where the last inequality follows from $r>2+\sqrt{2(2\alpha +\alpha c_2+c_1)\over 1-\alpha}$. The proof of Theorem 2 is complete.
\bigskip
\noindent{\it Proof of Theorem 3.} Since $\Pi_n(A_{r_n\varepsilon_n})\leq \Pi_n\bigl({\cal G}_n\cap A_{r_n\varepsilon_n}\bigr)+\Pi_n\bigl(A_{\varepsilon_n}\setminus {\cal G}_n\bigr),$ it suffices that the terms on the right hand side both tend to zero in probability. Given $\delta>0$, the proof of Lemma 3 implies that
$$F_0^\infty\bigl\{\Pi_n\bigl(A_{\varepsilon_n}\setminus {\cal G}_n\bigr)\geq \delta\bigr\}\leq {e^{2n\varepsilon_n^2}\Pi(A_{\varepsilon_n}\setminus {\cal G}_n)\over \delta\  \Pi(B_{\varepsilon_n^2})}+{\rm o}(1),$$
which by condition (1) tends to zero as $n\to\infty$.
Assume that $[r_n]$ stands for the largest integer less than or equal to $r_n$ and assume that  $D_j=\bigl\{f\in {\cal G}_n:\,j\varepsilon_n\leq H(f_0,f)< 2j\varepsilon_n\bigr\}$ ( Indeed, $D_j$ is an empty set for $j>\sqrt{2}/\varepsilon_n$ since the Hellinger distance cannot exceed $\sqrt{2}$ ).
Then we have
$$F_0^\infty\bigl\{\Pi_n\bigl({\cal G}_n\cap A_{r_n\varepsilon_n}\bigr)\geq \delta\bigr\}\leq F_0^\infty\Bigl\{\,\sum\limits_{j=[r_n]}^\infty\Pi_n(D_j)\geq \delta\,\Bigr\}$$ $$\leq F_0^\infty\Bigl\{\,\sum\limits_{j=[r_n]}^\infty\Pi_n(D_j)^\alpha \geq \delta\,\Bigr\}
\leq {1\over \delta}\,E\,\sum\limits_{j=[r_n]}^\infty\Pi_n(D_j)^\alpha
\leq {1\over \delta}\,\sum\limits_{j=[r_n]}^\infty E\, \Pi_n(D_j)^\alpha.$$
Take a partition $\bigcup_{i=1}^{N_j}D_{ji}$ for each $D_j$ such that $D_{ji}\subset \{f:\,H(f_{ji},\,f)<{j\varepsilon_n\over 3}\}$ for some $f_{ji}$ in $L_\mu$ and
$$\sum_{i=1}^{N_j} \Pi(D_{ji})^\alpha\leq 2\,\exp{\Bigl( J\bigl({j\varepsilon_n\over 3},D_j,\alpha\bigr) \Bigr)}\leq 2e^{c_1j^2n\varepsilon_n^2}\,  \Pi(B_{\varepsilon_n^2})^{\alpha},$$
where the last inequality follows from condition (2). Using the same argument as the proof of Theorem 1, one can get $H(f_{kD_{ji}},f_0)\geq {j\varepsilon_n/3}$ and hence we have with probability tending to 1
$$F_0^\infty\bigl\{\Pi_n\bigl({\cal G}_n\cap A_{r_n\varepsilon_n}\bigr)\geq \delta\bigr\} \leq {1\over \delta}\, \sum\limits_{j=[r_n]}^\infty E\,\Bigl(\ \sum_{i=1}^{N_j} \Pi_n(D_{ji})\, \Bigr)^\alpha$$
$$\leq {1\over \delta}\, \sum\limits_{j=[r_n]}^\infty  \sum_{i=1}^{N_j} E\,\Pi_n(D_{ji})^\alpha \leq  {e^{2\alpha n\varepsilon_n^2}\over  \delta\,\Pi(B_{\varepsilon_n^2})^\alpha}\,\sum\limits_{j=[r_n]}^\infty  \sum_{i=1}^{N_j} \Pi(D_{ji})^\alpha \,e^{{1\over 18}\,j^2n\varepsilon_n^2(\alpha-1)}$$
$$\leq {2\over \delta}\,\sum\limits_{j=[r_n]}^\infty   e^{n\varepsilon_n^2(2\alpha+c_1j^2+{1\over 18}\,j^2(\alpha-1))}
\leq {2\over \delta}\,\sum\limits_{j=[r_n]}^\infty   {1\over -n\,\varepsilon_n^2\,\bigl(2\alpha+c_1\,j^2+{1\over 18}\,j^2\,(\alpha-1)\bigr)}$$
$$\leq {3\over \delta \,n\,\varepsilon_n^2}\,\sum\limits_{j=[r_n]}^\infty   {1\over -c_1\,j^2+{1\over 18}\,j^2\,(1-\alpha)}\leq {3\over \delta \,n\,\varepsilon_n^2\,\bigl(-c_1+{1\over 18}\,(1-\alpha)\bigr)}\,\sum\limits_{j=[r_n]}^\infty   \Bigl({1\over j-1}-{1\over j}\Bigr)$$
$$={3\over \delta \,n\,\varepsilon_n^2\,\bigl(-c_1+{1\over 18}\,(1-\alpha)\bigr)([r_n]-1)}
$$
which tends to zero as $n\to\infty$, since $n\,\varepsilon_n^2$ are uniformly bounded away from zero and $r_n\to\infty$, where the second inequality follows from $\alpha\in (0,1)$, the third from the proof of Theorem 1, the fifth from the elementary inequality $e^{-x}<{1\over x}$ for $x>0$ and some of the inequalities only hold for all large $n$. Thus, we have proved that $\Pi_n\bigl({\cal G}_n\cap A_{r_n\varepsilon_n}\bigr)$ converges to zero in probability. The proof of Theorem 3 is complete.
\bigskip
\noindent{\it Proof of Theorem 4.} From $n\,\varepsilon_n^2\geq c_0\,\log n$ and $c_0\,c_2>1$ it turns out that $e^{-n\,\varepsilon_n^2\,c_2}\leq 1/n^{c_0\,c_2}$ and hence, by the first Borel-Cantelli Lemma and Lemma 1, we obtain that
$$\int_{\F} R_n(f)\, \Pi(df)\geq e^{-n\,\varepsilon_n^2\,(3+2c_2)}\ \Pi\bigl( W_{\varepsilon_n})$$
almost surely for all large $n$.  Then, following the proofs of Lemma 3 and Theorem 3, one can get that for any $\delta>0$ and $r>1$,
$$F_0^\infty\bigl\{\Pi_n\bigl( A_{r\varepsilon_n}\bigr)\geq \delta\bigr\} \leq F_0^\infty\bigl\{\Pi_n\bigl( A_{\varepsilon_n}\setminus {\cal G}_n\bigr)\geq {\delta\over 2}\bigr\}+F_0^\infty\bigl\{\Pi_n\bigl({\cal G}_n\cap A_{r\varepsilon_n}\bigr)\geq {\delta\over 2}\bigr\}$$
$$\leq {2\, e^{n\,\varepsilon_n^2\,(3+2c_2)}\,\Pi(A_{\varepsilon_n}\setminus {\cal G}_n\bigr)\over \delta\,\Pi(W_{\varepsilon_n})}+{4\over \delta}\,\sum\limits_{j=[r]}^\infty   e^{n\varepsilon_n^2((3+2c_2)\alpha+c_1j^2+{1\over 18}\,j^2(\alpha-1))}:=a_n+b_n.$$
Condition (1) yields that $\sum_{n=1}^\infty a_n<\infty$. On the other hand, since $c_1< {1-\alpha\over 18}$, we have that for all $n\geq 2$ and for all $r$ so large that $(3+2c_2)\alpha+(c_1+{\alpha-1\over 18})\,(r^2-1)\leq -{2\over c_0}$,
$$b_n\leq {4\, e^{n\varepsilon_n^2(3+2c_2)\alpha}\over \delta}\,\sum\limits_{j=[r]^2}^\infty   e^{n\varepsilon_n^2(c_1+{\alpha-1\over 18})\,j}={4\, e^{n\varepsilon_n^2\big((3+2c_2)\alpha+(c_1+{\alpha-1\over 18})\,[r]^2\big)}\over \delta\,\big(1-e^{n\varepsilon_n^2(c_1+{\alpha-1\over 18})}\big)}$$
$$\leq {4\, n^{c_0\big((3+2c_2)\alpha+(c_1+{\alpha-1\over 18})\,[r]^2\big)}\over \delta\,\big(1-n^{c_0(c_1+{\alpha-1\over 18})}\big)}\leq {4\, n^{c_0\big((3+2c_2)\alpha+(c_1+{\alpha-1\over 18})\,(r^2-1)\big)}\over \delta\,\big(1-2^{c_0(c_1+{\alpha-1\over 18})}\big)}\leq {4\, n^{-2}\over \delta\,\big(1-2^{c_0(c_1+{\alpha-1\over 18})}\big)}.$$
and hence $\sum_{n=1}^\infty b_n<\infty$. Thus, by the first Borel-Cantelli Lemma  we obtain that $\Pi_n\bigl( A_{r\varepsilon_n}\bigr)\leq \delta$
almost surely for all large $n$, which concludes the proof of Theorem 4.
\bigskip \bigskip
\bigskip \centerline{ REFERENCES } \bigskip

\noindent B{\eightrm ARRON}, A., S{\eightrm CHERVISH}, M. and W{\eightrm ASSERMAN}, L. (1999). {\eightrm The consistency of posterior

\qquad distributions in nonparametric problems. Ann. Statist. {\bf 27}, 536-561. }

\noindent G{\eightrm HOSAL}, S.  (2001). {\eightrm Convergence rates for density estimation with Bernstein polynomials. Ann.

\qquad  Statist. {\bf 29}, 1264-1280. }

\noindent G{\eightrm HOSAL}, S., G{\eightrm HOSH}, J. K.  and R{\eightrm AMAMOORTHI}, R. V. (1997). {\eightrm Non-informative priors via

\qquad sieves  and packing numbers. In Advances in Statistical Decision Theory and Applications (S.

\qquad Panchapakeshan and N.Balakrishnan eds.) 129-140. Birkh\"auser, Boston.}

\noindent G{\eightrm HOSAL}, S., G{\eightrm HOSH}, J. K.  and V{\eightrm AN DER }V{\eightrm AART}, A. W. (2000). {\eightrm Convergence rates of posterior

\qquad distributions. Ann. Statist. {\bf 28}, 500-531. }

\noindent G{\eightrm HOSAL}, S.  and V{\eightrm AN DER }V{\eightrm AART}, A. W. (2001). {\eightrm Entropies and rates of convergence for maximum

\qquad likelihood and Bayes estimation for mixtures of normal densities.  Ann. Statist. {\bf 29}, 1233-1263. }

\noindent G{\eightrm HOSAL}, S.  and V{\eightrm AN DER }V{\eightrm AART}, A. W. (2007a). {\eightrm Convergence rates of posterior distributions

\qquad for noniid observations.  Ann. Statist. {\bf 35}, 192-223. }

\noindent G{\eightrm HOSAL}, S.  and V{\eightrm AN DER }V{\eightrm AART}, A. W. (2007b). {\eightrm Posterior convergence rates of Dirichlet

\qquad  mixtures at smooth densities. Ann. Statist. {\bf 35}, 697-723. }

\noindent K{\eightrm OLMOGOROV}, A. N.  and T{\eightrm IHOMIROV}, V. M. (1959). {\eightrm $\varepsilon$-entropy and $\varepsilon$-capacity of sets in

\qquad function spaces. Uspekhi Mat. Nauk {\bf 14}, 3-86 [in Russian; English transl. Amer. Math. Soc.

\qquad  Transl. Ser. 2, {\bf 17}, 277-364 (1961)].}

\noindent L{\eightrm IJOI}, A., P{\eightrm RUNSTER}, I. and W{\eightrm ALKER}, S.  (2007). {\eightrm On convergence rates for nonparametric

\qquad  posterior distributions. Aust. N. Z. J. Stat. {\bf 49}$\,$(3), 209-219. }

\noindent P{\eightrm ETRONE}, S.  (1999). {\eightrm Random Bernstein polynomials. Scand. J. Statist. {\bf 26}, 373-393.}

\noindent P{\eightrm ETRONE}, S. and W{\eightrm ASSERMAN}, L.  (2002). {\eightrm Consistency of Bernstein polynomial posteriors. J. R.

 \qquad Stat. Soc. Ser. B Stat. Methodol.  {\bf 64}, 79-100.}

 \noindent S{\eightrm CHWARTZ}, L. (1965). {\eightrm On Bayes procedures Z. Wahr. verw. Geb. {\bf 4}, 10-26. }

 \noindent S{\eightrm CRICCIOLO}, C.  (2006). {\eightrm Convergence rates for Bayesian density estimation of infinite-dimensional

\qquad  exponential families. Ann. Statist. {\bf 34}, 2897-2920. }

\noindent S{\eightrm HEN}, X. and W{\eightrm ASSERMAN}, L.  (2001). {\eightrm Rates of convergence of posterior distributions. Ann.

\qquad  Statist. {\bf 29}, 687-714. }

\noindent S{\eightrm TONE}, C. J. (1990). {\eightrm Lerge-sample inference for log-spline models. Ann. Statist. {\bf 18}, 717-741. }

\noindent  W{\eightrm ALKER}, S.  (2004). {\eightrm New approaches to Bayesian consistency. Ann. Statist. {\bf 32}, 2028-2043. }

\noindent  W{\eightrm ALKER}, S., L{\eightrm IJOI}, A. and P{\eightrm RUNSTER}, I. (2007). {\eightrm On rates of convergence for posterior

\qquad distributions in infinite-dimensional models. Ann. Statist. {\bf 35}, 738-746. }

\noindent X{\eightrm ING}, Y. and R{\eightrm ANNEBY}, B.  (2008). {\eightrm Sufficient conditions for Bayesian consistency. (Preprint). }
\bigskip\bigskip
\noindent Yang Xing

\noindent Centre of Biostochastics

\noindent Swedish University of Agricultural Sciences

\noindent  SE-901 83, Ume\aa

\noindent  Sweden

\noindent E-mail address:\enskip yang.xing@sekon.slu.se

\end